\documentclass[preprint,12pt,authoryear]{elsarticle}

\usepackage{amssymb}
\usepackage{float}

\usepackage{algorithm}
\usepackage{algorithmic}

\begin{document}

\begin{frontmatter}



\title{ An approximate operator-based learning method for the numerical solution of stochastic differential equations}


\author[1]{Jingyuan Li}
\affiliation[1]{organization={School of Mathematics and Statistics},
            addressline={Wuhan University}, 
            city={Wuhan},
            country={China}}

\author[1]{Wei Liu \corref{cor1}}

\cortext[cor1]{Corresponding author}

\begin{abstract}
Stochastic differential equation (SDE in short) solvers find numerous applications across various fields. However, in practical simulations, we usually resort to using Itô-Taylor series-based methods like the Euler-Maruyama method. These methods often suffer from the limitation of fixed time scales and recalculations for different Brownian motions, which lead to computational inefficiency, especially in generative and sampling models. To address these issues, we propose a novel approach: learning a mapping between the solution of SDE and corresponding Brownian motion. This mapping exhibits versatility across different scales and requires minimal paths for training. Specifically, we employ the DeepONet method to learn a nonlinear mapping. And we also assess the efficiency of this method through simulations conducted at varying time scales. Additionally, we evaluate its generalization performance to verify its good versatility in different scenarios.

\end{abstract}




\end{frontmatter}


\section{Introduction}
Stochastic differential equations (SDEs) are one of the most important representations of dynamical systems. They are notable for the ability to include a deterministic component of the system and a stochastic one to represent random unknown factors and can be used to model phenomena from a variety of disciplines ranging from finance to hydrology, from rotational diffusion in granular media to climatology \citep{oksendal2013stochastic,mao2007stochastic}.

Nowadays SDEs has increasing importance in both theory and applications. Especially, numerical solvers for SDEs have become important recently due to the extensive use of diffusion models \citep{Yang2022DiffusionMA}, generative models \citep{creswell2018generative}, SDE-based MCMC algorithms \citep{girolami2011riemann,durmus2019analysis}, etc. For example, the convergence rate to the invariant distribution in Langevin MCMC sampling often depend on the numerical discretization scheme.  However, in practical applications we seldom use continuous models directly because the traditional numerical solvers often require very small discretization errors and a large amount of samples to obtain the corresponding statistical estimates, which lead to huge computational cost. Indeed due to the exponential explosion of the It\^o-Taylor series-based methods, we often have to trade off between the size of iteration step and the time length of the algorithm run. And for the simulation of more different paths we can not directly use the simulated results we already obtained, instead we need to iterate the algorithm from the beginning again.

In this paper we only consider the SDEs drive by Brownian motions. It is known that by Watanabe Theorem \citep{ikeda2014stochastic}, under certain assumptions on the coefficients of the SDE, the solution of SDE can be represented as a mapping (functional) on the initial value and the corresponding Brownian motion. Explicitly, for a well-defined SDE with strong solutions
\begin{equation}\label{SDE} dX_t = a(t,X_t)dt + b(t,X_t)dB_t \quad X_0=x_0 \end{equation}
where $a:\left[0,T\right] \times  {\bf R}^d \to {\bf R}^d$ is a vector-valued function, $b: \left[0,T\right] \times {\bf R}^d\to {\bf R}^d\times {\bf R}^m$ is a matrix-valued function, and $B$ is some $m$-dimensional Brownian motion. When the coefficients $b$ and $\sigma$ satisfy some regularity conditions, the solution can be represented as the following mapping (or functional) $\mathbf{F}$ on the Brownian motion and the initial value $X_0$:
\begin{equation}\label{Map}X_t = \mathbf{F}(X_0,B)(t).\end{equation}
 By the classical results of stochastic analysis \citep{Kallenberg2021FoundationsOM}, the existence of this representation often does not require very strong assumptions such as the Lipschitz condition, but only the existence of weak solutions to the SDE and the path uniqueness property. This representation can  be used to analyze the property of the solution. However, it is not easy to use this representation to design an algorithms for solving SDE.
Recently the DeepONet \citep{lu2019DeepONet} provided a way to learn the nonlinear mapping by only a portion of samples and then use it to make prediction.

\textbf{Our contributions:} we propose a general framework to learn the mapping $\mathbf{F}$ in (\ref{Map}) between the solution of SDE (\ref{SDE}) and associated  Brownian motion $B$. This approach has great potential applications since it can reduce significantly the computation cost and improve the efficiency (compared to the existing SDE solvers). We present the main advantages explicitly in the following:

\begin{enumerate} 
\item By using DeepONet to learn the mapping $\mathbf{F}$, we can get better generalisation results i.e. better results with less samples.

\item  This approach is meshless and suitable for {\bf varying time scales/ multiscale}. This means that the mapping $\mathbf{F}$ we learn by using any specification of mesh has generalisation abilities on other scales of mesh as well, since DeepONet is meshless. This model has good predictive ability. We consider Brownian motions with time scales $\left [0,h,2h.... ,(M-1)h\right ]$, where $T = (M-1)h$, and train the model at a larger step size $h$. This model can be used for prediction at smaller step size $h$. The converse is also possible. For the same step size, we can use one part of the Brownian motion paths for training and then the other part of the Brownian motion paths as prediction.

\item This approach can greatly improve the efficiency of high-precision SDE solvers. We only need a small number of high-precision training samples to learn the mapping, and then use this mapping to generate directly more paths of the solution of SDE.  And we do not even need to simulate the entire path of Brownian motion, we only need to simulate the value of Brownian motion at the time points we are interested in.

\item This approach can be applied to a wide range of SDEs. For instance we consider the nonlinear McKean-Vlasov SDE which is difficult to compute or simulate its solution (due to the high nonlinearity and distribution dependence). Our approach achieves a complete reconstruction of the McKean-Vlasov SDE, and thus successfully avoids to using the particle system to approximate McKean-Vlasov SDE. By doing so we reduce significantly the computation cost, which makes the McKean-Vlasov SDE more applicable in practice. Moreover, besides the SDEs driven by Brownian motions, this approach can be also applied SDEs driven by  Poisson Random Measure \citep{Zhao2014YamadaWatanabeTF} and SPDEs \citep{Schmuland2008YamadaWatanabeTF,Tappe2013TheYT} etc. This greatly expands the use of this methodology. 

\item  This approach can somewhat unify continuous stochastic differential equation models and normalized flow methods in the field of sampling, and thus may contribute to a better understanding of normalized flow methods.

\end{enumerate}
We validate the effectiveness of this approach through doing some experiments. This method may greatly promote the development of related fields, such as fast sample generation, fast sampling, etc.
\section{Related Work}
A large class of generative models based on score \citep{dockhorn2021score,song2020score} is the class of generative models that often require matching scores and then sampling the target distribution using an SDE solver. 
Normalized flow \citep{zhang2021diffusion,rezende2015variational} is another generative model that takes a prior distribution $z$ and turns it into the target distribution $x$ by a transformation $f$.
In the field of SDE solvers, the commonly used methods are Itô-Taylor series-based methods \citep{kloeden1992stochastic}. 
The strengths of the DeepONet approach are its powerful modelling capabilities efficient training and inference capabilities, and its wide range of applicability, making it an effective tool for dealing with complex problems. For the learning of operators for ODEs with control \citep{Tan2022EnhancedDF}. There are aspects of prediction for general stochastic systems \citep{Garg2022AssessmentOD,zhang2022multiauto}. DeepONet for transfer learning \citep{Goswami2022DeepTL}. Solving the inverse problem \citep{Zhang2022MultiAutoDeepONetAM}.
emulation of parametric differential equations \citep{Wang2021LearningTS}. In prediction of multi-scale systems \citep{Lin2020OperatorLF}. Error estimates for DeepOnets \citep{Lanthaler2021ErrorEF}.

\section{Theory and Algorithms}
\subsection{The solutions of SDE}
In the following we first present some known results about the solution of SDE.

 \subsubsection{Strong solution:}Let $\{B_t\}_{t\ge 0}$ be a standard Brownian motion on a probability space $(\Omega,\mathcal{F},P)$ and let $\{\mathcal{F}_t\}_{t\ge 0}$ be the completion of the minimal filtration by null sets. A {\it strong solution} $\{X_t\}_{t\ge0}$ of the SDE (\ref{SDE}) with initial value $X_0$ is a continuous and adapted process such that for all $t\ge0$,
\begin{equation}\label{sol}X_t=X_0+\int_0^t a(s,X_s)ds+\int_o^t b(s,X_s)dB_s\ a.s. \end{equation}

\subsubsection{Weak solution:} A {\it weak solution} of the SDE (\ref{SDE}) with initial value $X_0$ is a continuous stochastic process $X_t$ defined on {\it some} probability space $(\Omega,\mathcal{F},P)$  such that for some  Brownian motion and some admissible filtration the process $X_t$ is adapted and satisfies the stochastic integral equation (\ref{sol}).

What distinguishes a strong solution from a weak solution is the requirement that it be adapted
to the completion of the minimal filtration. This makes each random variable $X_t$ a measurable
function of the path $\{B_s\}_{s\le t}$ of the driving Brownian motion. It turns out that for
certain coefficient functions $a$ and $b$, solutions to the stochastic integral equation equation (\ref{sol}) may
exist for some Brownian motions and some admissible filtrations but not for others.
 
\subsubsection{Another form of weak solution:}
The distribution $\bf{P}$ of the solution $X_t$ is solved by the local martingale problem:
\[M_t^f = f(X_t)-f(X_0)-\int_{0}^{t} A_s f(X_s )ds \quad f \in C^\infty
\] 
\[A_s f(x)= \frac{1}{2}\left [{b(s,x)b(s,x)^T}\right ]^{ij}f''_{ij}(x_s) + a^i(s,x)f'_j(x_s)\]
For a given $(a,b)$, $M_t^f$ is local martingale.

\subsubsection{Pathwise uniqueness of the solution:}
Pathwise uniqueness is said to hold if $X$ and $Y$ are two solutions of SDE (\ref{SDE}) defined on the same probability space and $X_0 = Y_0$ a.s., then $X_t = Y_t$ a.s., $\forall t \geq 0$.

\subsubsection{Yamada-Watanabe Theorem :}\citep{Kallenberg2021FoundationsOM}
Let the functions  $a$, $b$  be progressive and such that weak existence and pathwise uniqueness hold for solutions to equation  ($a$, $b$)  starting at fixed points. Then

(i) strong existence and uniqueness in law hold for any initial distribution.

(ii) there exists a measurable, universally adapted function $F$ , such that every solution  ($X$, $B$)  to equation  ($a$, $b$)  satisfies  $X=F\left(X_{0}, B\right)$  a.s.

\subsubsection{Classical regularity conditions:} The following sufficient condition is often used in It\^o-Taylor series-based methods. Note that this condition is a lit stronger and can be weakened (locally Lipschitz).

Lipschitz condition: there exists some constant $L>0$ such that
\[|a(t, x)-a(t, y)| \leq L|x-y|\]
\[|b(t, x)-b(t, y)| \leq L|x-y|\]
for all  $t \in\left[0, T\right]$ and  $x, y \in {\bf R}^d$.

Linear growth bound: there exists some constant $C>0$ such that
\[|a(t, x)|^{2} \leq C (1+|x|^{2})\]
\[|b(t, x)|^{2} \leq C (1+|x|^{2})\]
for all  $t \in\left[0, T\right]$ and  $x, y \in {\bf R}^d$.

\subsubsection{SDE solver:}
 It\^o-Taylor series-based methods are based on the It\^o-Taylor expansion of stochastic integral, and the convergence order is usually related to the approximation order of the stochastic integral. The most commonly used method is the Euler-Maruyama discretization, which admits the following scheme:
\[X_{t+h} =X_t + a(t,X_t)h+b(t,X_t)(B_{t+h}-B_t)\]
The error can be controlled by the following result under some regularity conditions:
\[E\left [ \sup_{t\in[0,T]}|X_t - \hat X_t|^2 \right ] \leq C_T h\]
where $\hat X_t$ is the solution obtained by the SDE solver.

Since the constant $C_T$ usually depends exponentially on the time period $T$, these numerical approaches exhibit poor accuracy for large $T$. However in some practical cases we need to make $T$ very large, e.g. approximating the stationary distribution of SDE (Langevin Monte-Carlo). To reduce this error we will take smaller step size $h$ and simulate more paths which cost too much computation resource.

However, if we can use a smaller $h$ to obtain solutions for a portion of the simulation, although this will take longer, we can use these solutions to get a mapping of $B_t$ with respect to $X_t$, and then simulate a portion of the Brownian motion, which can then be predicted directly using the model can greatly improve the efficiency of the SDE solver.

\subsection{Approximate the mapping $\bf{F}$ by DeepONet}
DeepONet is a network that can approximate continuous nonlinear operators. The architecture of the network mainly contains two parts: one or several branch networks and trunk network. It can be written in the following form :
\begin{figure*}[t]
\centering
\includegraphics[width=0.8\textwidth]{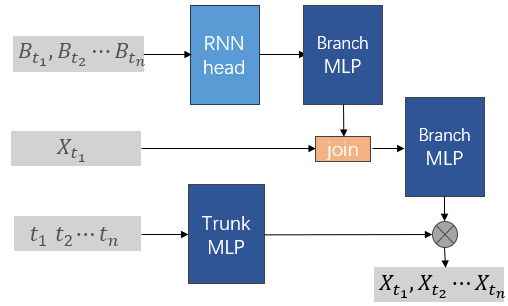}
\caption{DeepONet's architecture approximates mapping $F$}
\label{fig1}
\end{figure*}

\begin{small}
\[
G(u)(y) \approx \sum_{k=1}^{p} \underbrace{\sum_{i=1}^{n} c_{i}^{k} \sigma\left(\sum_{j=1}^{m} \xi_{i j}^{k} u\left(x_{j}\right)+\theta_{i}^{k}\right)}_{\textit {Branch }} \underbrace{\sigma\left(w_{k} \cdot y+\zeta_{k}\right)}_{\textit {Trunk }}
\]
\end{small}
where $c_{i}^{k},\xi_{i j}^{k},\theta_{i}^{k}w_{k},\zeta_{k}$ are parameters that can be trained. The $\sigma$ is an activation function.  However, in the setup of SDE, this may not express the mapping ${\bf F}$ since $X_t$ may often depend on the path of the Brownian motion before $t$ i.e. $\{B_s\}_{0\le s\le t}$, which leads to some \textbf{memory property}. Therefore we need to learn this memory property by using an RNN head. We first map the samples into a hidden space and then use a branch network to realize that.
We set $u$ in the branch network to $MLP(RNN(B),X_0)$ and $y$ in trunk networks to $t$.
Of course the RNN head is not the only choice. We can also consider using other Recurrent neural network network structures, such as LSTM \citep{Hochreiter1997LongSM}.

Unlike the usual DeepONet, we need to add a term about the initial value dependence in the architecture of the network, which is made by splicing the values directly after the $B_t$ mapping, then the vector is mapped with the MLP after the splicing is completed. For the calibration of the initial value, we use the mean square error of all the points on the track, including the $t=t_1$ moment, as the loss function. For details, we can refer to Figure \ref{fig1}.

We use the mean square error as the loss function: for a given number of sample paths $(B_t^i,X_t^i)$, $1 \leq i \leq N$, $0 \leq t_k \leq T=(M-1)h$, $t_{k+1} - t_k = h$ and $X_t^i$ is simulated numerically from the SDE solver,
\[
\bf{F}_{\theta^*} = \arg\min_{\bf{F}_{\theta} }\sum_{i=1}^{N} \sum_{k=1}^{M} ({\bf F}_{\theta}(X_0^i,B^i)(t_k)-X_{t_k}^i)^2
\]
where $\bf{F}_\theta$ is DeepONet. We have the Algorithm \ref{alg:algorithm} which is presented later.

\subsection{Some examples}
\subsubsection{Some examples of SDEs of the It\^o type:}

Consider the following SDE:
\[dX_t =aX_t dt + bX_tdB_t\]
\[X_t = X_0\exp((a-\frac{1}{2}b^2)t+bB_t.\]
The solution of the SDE is the Geometric Brownian motion. In this case, $X_t$ depends only on $B_t$ which is the value of the Brownian motion at current moment $t$, thus it will be easy to learn ${\bf F}$ as a multivariate functional.

For another class of cases it can be shown that not only does it depend on the $B_t$ at the current moment.
\[dX_t =-aX_t dt + bdB_t\]
\[X_t = \exp(-at)(X_0+b \int_{0}^{t}\exp(as)dB_s )\]
Unlike the case of geometric Brownian motion, here $X_t$ depends on the path of the Brownian motion on the time interval $[0,t]$.

In the area of sampling, we have the Langevin MCMC algorithm with the target density function as following:
\[p(x)= C\exp(-V(x))\]
 where $C$ is the normalized constant that does not need to be known in advance. This target distribution is just 
 the stationary distribution of the following SDE:
\[dX_t= \frac{1}{2}\nabla \log(p(x))dt + dB_t.\]
When $T$ is sufficiently large we can approximate the samples of the distribution $p(x)$.

\begin{algorithm}[tb]
\caption{SDE operator fitting algorithm}
\label{alg:algorithm}
\textbf{Input}: Brownian motion $B_t^i$ and Solution of SDE $X_t^i$ where $1 \leq i \leq N$ and $0 \leq t_k\leq T=(M-1)h$. $M$ is the number of Euler steps and $h$ is the step length\\
\textbf{Parameter}: Parameters $\theta$ for DeepONet network ${\bf F}_{\theta}$\\
\textbf{Output}: ${\bf F}_{\theta}(B^i,X_0)(t)$
\begin{algorithmic}[1] 
\STATE Calculate the loss 
\[
\textit {loss }=\frac{1}{MN}\sum_{i=1}^{N} \sum_{k=1}^{M} ({\bf F}_{\theta}(X_0^i,B^i)(t_k)-X_{t_k}^i)^2
\]
\STATE Update network parameters $\theta$ using optimizer
\IF {$\textit {loss } \leq \textit {Threshold}$}
\STATE End of training
\ENDIF
\STATE \textbf{return} DeepONet network ${\bf F}_{\theta}(X_0,B^i)(t)$
\end{algorithmic}
\end{algorithm}

\subsubsection{McKean-Vlasov SDE}
We consider a nonlinear SDE with high numerical computational complexity. McKean-Vlasov SDE is a class of stochastic differential equations describing the interactions between particles and affected by random perturbations. This class of SDEs is characterised by the fact that the evolution of each particle depends not only on its own state and external stochastic perturbations, but is also influenced by the average behaviour of the whole particle system. Specifically, the evolution of each particle includes in its dynamics the influence of the overall particle density distribution. Explicitly it has the following form:
\[
d X_{t}=a\left(t, X_{t}, \mu_{t}\right) d t+b\left(t, X_{t}, \mu_{t}\right) dB_{t}
\]
where $\mu_{t}$ is the distribution of $X_t$. The related Yamada-Watanabe Theorem has \citep{grube2023strong}
McKean-Vlasov SDE is the mean field limit of the following particles system:
\[
d X^n_{t}=a\left(t, X^n_{t}, \hat \mu_{t}\right) d t+b\left(t, X^n_{t},\hat \mu_{t}\right) d B^n_{t}
\]
where $\hat{\mu}^{N}_t:=\frac{1}{N} \sum_{n=1}^{N}\delta_{X^n_t}$ is empirical distribution of $X_t$ and $\delta_x$ denotes the Dirac mass at point $x$.

We usually use the standard SDE solver for particle-particle systems, and then approximate the McKean-Vlasov SDE (due to propagation of chaos). Some of the existing SDE solver \citep{Antonelli2002RateOC,reisinger2023convergence,Liu2022ParticleMA,Bao2021FirstorderCO}. It is often very difficult to solve this equation. This is mainly due to the inefficiency of numerical computation for several reasons:

\textbf{Distribution-dependent terms :} The evolution of such SDEs also depends on the distribution of the overall particle density, not just the state of individual particles. This requires approximation of the distribution-dependent terms in the numerical calculations, which increases the computational complexity. Since it involves the estimation of the overall particle density, empirical measures at each time step is used to approximate the whole overall particle density, which leads to high computational costs.

\textbf{Computational cost of SDE numerical solvers:} To solve the McKean-Vlasov SDE, we also need to use standard SDE solvers, such as the Euler-Maruyama method. Additionally, it is necessary to incorporate techniques for approximating empirical measures within the SDE solver. Hence, to ensure adequate precision, it becomes necessary to take small time intervals. This increases the computational effort by requiring a large number of time steps to be performed in a limited amount of time. In particular, when calculating stationary distributions, it is often necessary to perform a large number of steps and a large number of samples.

Our approach can avoid the use of empirical measures to approximate distribution-dependent terms by using a limited number of high-quality solution trajectories to learn to map $\bf{F}$ with respect to $X_t$ and $(X_0,B_t)$. Hence, to obtain the solution $X_t$ once more, we can directly implement the original approximation of the McKean-Vlasov SDE and avoid the use of particle systems. At the same time, the significant reduction in computational requirements greatly facilitates the fast sampling of this particular type of SDE.

\section{Experiments}
\subsection{Generalization performance of the model}
Our first example considers the following SDE to express the generalization performance of the model.
\[dX_t =-aX_t dt + bdB_t\]
The SDE solution is 
\[X_t = \exp(-at)\left(X_0+b \int_{0}^{t}\exp(as)dB_s \right)\]

\begin{figure}[t]
\centering
\includegraphics[width=0.90\columnwidth]{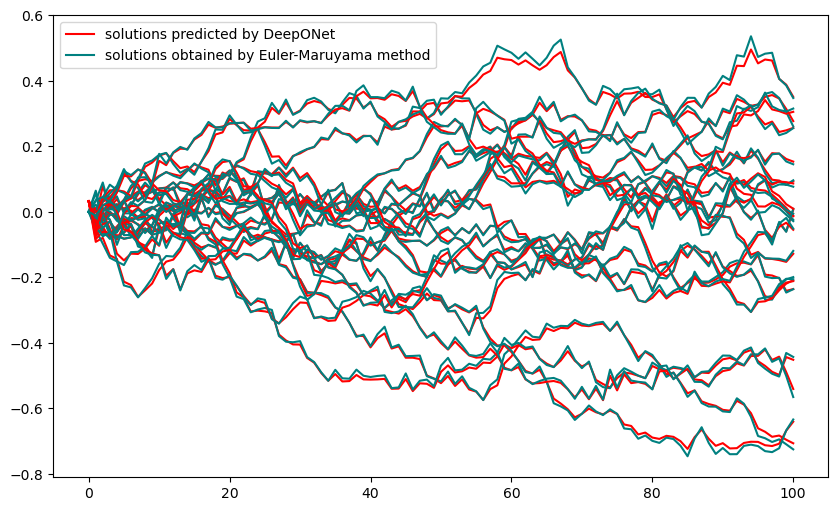}
\includegraphics[width=0.90\columnwidth]{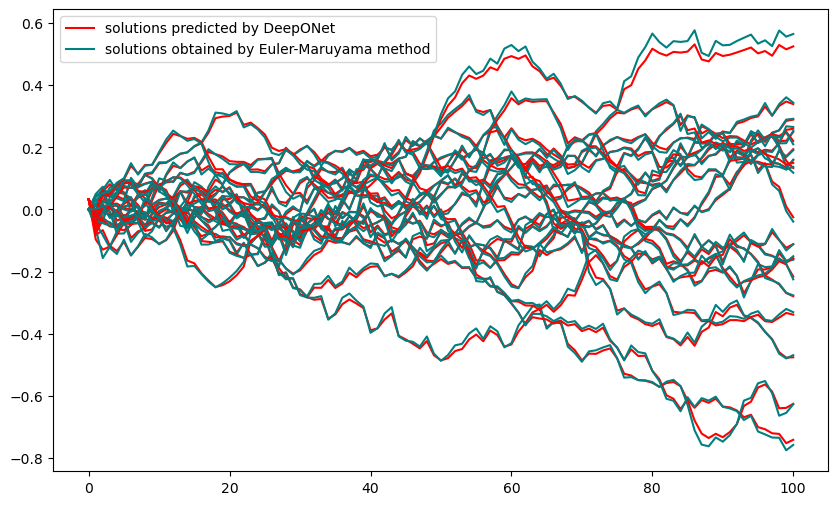} 
\caption{Paths obtained using Euler-Maruyama method and DeepONet on the training and test sets.The training set was trained using 20 paths. The top is the training set and the bottom is the test set}
\label{fig2}
\end{figure}

\begin{figure}[t]
\centering
\includegraphics[width=0.95\columnwidth]{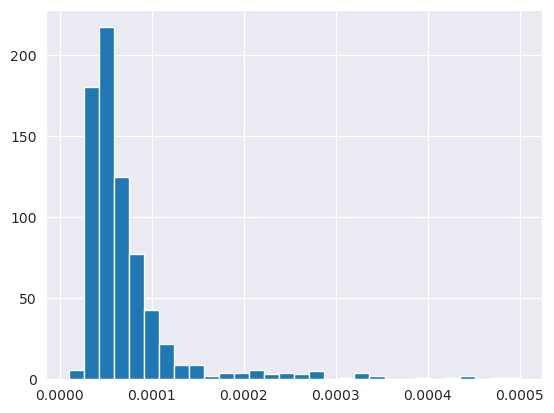}
\caption{Distribution of prediction errors}
\label{fig3}
\end{figure}

We used 20 paths $(B_t,X_t)$ for training and then used the trained model to predict. Each of these generated tracks of Brownian motion corresponds to a random number seed. The seeds for training and prediction do not overlap. We get direct output from the model. And we predict the solution $X_t$ for 800 paths and compute the Mean Squared Error (MSE) with respect to the true solution. we plot the error in a Figure \ref{fig3}, where the range of the Figure \ref{fig3} is $\left[1e-5,5e-4\right]$.
We presented the predictive performance of a subset of projected paths, along with the effectiveness of fitting to the training dataset in Figure \ref{fig2}. It is evident that we achieved an approximation to the SDE operator with a minimal sample size of $N=20$. Moreover, we successfully attained reduced errors for larger sample sizes, indicating improved accuracy in our results.


Subsequently, for this particular example, we extend our analysis to encompass predictive capabilities across multiple scales. Specifically, we apply the previously obtained model to predict solutions of SDE corresponding to Brownian motions at various scales. We consider the following scales: $\frac{h}{1000},\frac{h}{100},\frac{h}{10},10h,100h,1000h$. Similarly, we performed predictions on 800 paths and calculated the Standardized Mean Squared Error (MSE) between the predicted and true values. This standardized metric scales both the predicted and true values to the range of 0 to 1 before computing the MSE, enabling a more comparable assessment of the prediction performance.
For all the errors, we calculated both the mean and standard deviation to quantify the discrepancies.
\begin{table}[t]
\centering
\begin{tabular}{ccc}
\hline
Scales & Mean & Standard deviation \\
\hline
$h$ & 1.335e-03 & 1.313e-03 \\
\hline
$\frac{h}{10}$  &7.831e-03 & 8.471e-03 \\
\hline
$\frac{h}{100}$  & 4.604e-02 & 2.743e-02 \\
\hline
$\frac{h}{1000}$  & 6.725e-02 & 2.479e-02 \\
\hline
$10h$  & 9.709e-03 & 1.133e-02 \\
\hline
$100h$  & 3.090e-01 & 8.988e-02 \\
\hline
$1000h$  &2.402e-01 & 1.311e-01 \\
\hline
\end{tabular}
\caption{Comparison of Errors at Different Scales}
\end{table}
\subsection{Fast sampling for the McKean-Vlasov SDE}

We consider the following SDE of the McKean-Vlasov type equation, which is a Burgers equation \citep{Bossy1997ASP}.
\[
d X_{t}=\int_{R}^{} H(X_t-y)d\mu_t(dy)dt+\sigma B_t \quad H(x) = I_{(x\geq 0)}
\]

We used a classical SDE solver, the Euler-Maruyama Particle (EMP) method
Briefly, given an initial value $X_0$, we obtain $X_T^{n, N}$ by the following iterative equation
\begin{small}
\[
X_{t_{m+1}}^{n, N} = X_{t_{m}}^{n, N} + h \, a\left(t_{m}, {X}_{t_{m}}^{n, N},\hat{\mu}_{t_{m}}^{N}\right)+ b\left(t_{m}, X_{t_{m}}^{n, N}, \hat{\mu}_{t_{m}}^{N}\right)\delta_{t}
\]
\end{small}
where $\delta_{t} =  \left(B_{t_{m+1}}^{n} - B_{t_{m}}^{n}\right)$,$\hat{\mu}_{t_{m}}^{N}:=\frac{1}{N} \sum_{n=1}^{N} \delta_{{X}_{t_{m}}^{n, N}}$

We set the step size $h = t_{m+1} -t_{m} = 0.01$, $\sigma= 1$, iterated for $M=31$ steps, and used $N = 10,000$ samples to obtain the ''true value'' of the solution.

\begin{figure}[t]
\centering
\includegraphics[width=0.95\columnwidth]{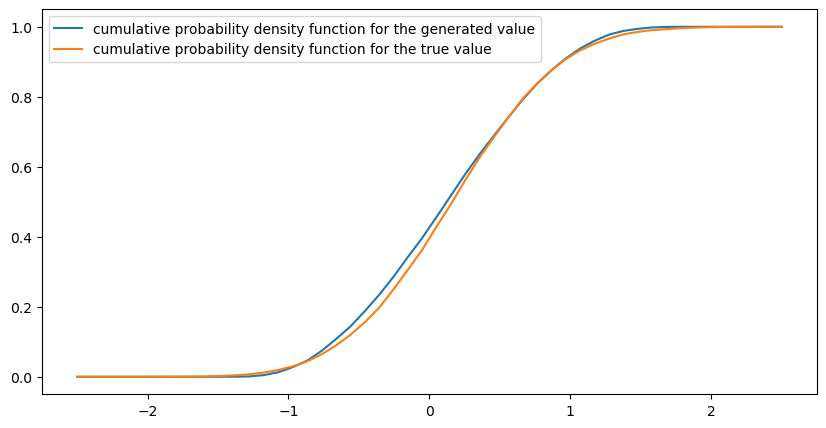}
\caption{cumulative probability density function curves}
\label{fig4}
\end{figure}

First, we pick 20 paths out of 10,000 samples to train $F$. Then using the learned $F$ to generate 10,000 paths of SDE solutions. We plot the cumulative probability density function curves in $X_T$ for the true values and the values of the generated solutions in Figure \ref{fig4}.

Further on we start to calculate the time needed for the two algorithms to generate the solution $X_T$, note that the time taken to simulate the Brownian motion is different for the two methods, in our method we can directly simulate the Brownian motion of a single point for the value $B_T$, whereas in EPM we need to simulate the whole path beforehand. All simulations are done on the same CPU (AMD EPYC 7B12). And we also recorded the time it took us to train ${\bf F}$ (T-time), and the time it took the model (O-time) to reason about it. And in the EPM approach, we calculated the time to simulate the entire Brownian motion path (B-time), and the time to compute it using the iterative formula (O-time). Specifically we summarise all the results in Table 2

\begin{table}[H]
\centering
\resizebox{1.0\textwidth}{!}{%
\begin{tabular}{cccccc}
\hline
Method & $N$ & $M$ & B-time & O-time & T-time \\
\hline
EMP & 100 & 31 & 5.6601e-03 & 2.4640e-02 &  $\times$    \\
DeepONet & 100 & 31 & 6.2943e-05 & 4.0356e-02 &  1.1468e+02 (20 paths)  \\ 

EMP & 1000 & 31 & 7.4580e-02& 2.9986e-01 &  $\times$    \\
DeepONet & 1000 & 31 & 3.5286e-05 & 4.0859e-01 &  1.1468e+02 (20 paths)  \\

EMP & 10000 & 31 & 4.8686e-01  & 9.1371e+00 &  $\times$    \\
DeepONet & 10000 & 31 & 3.7098e-04 & 6.8098e+00 &  1.1468e+02 (20 paths)  \\

EMP       & 100 & 51 &  1.1670e-02   & 4.3571e-02      &    $\times$    \\
DeepONet & 100 &  51 &   2.5272e-05  &  5.2540e-01   &  1.9591e+02 (20 paths)\\

EMP       & 1000 & 51 & 7.1276e-02    &  4.9766e-01     &    $\times$    \\
DeepONet & 1000 &  51 & 3.6955e-05    &   4.8736e-01    &  1.9591e+02 (20 paths)\\

EMP       & 10000 & 51 & 4.9503e-01    &   4.6647e+01   &    $\times$    \\
DeepONet & 10000 &  51 & 4.7207e-05    &   3.9963e+00   &  1.9591e+02 (20 paths) \\

EMP       & 100 & 101 &  1.1119e-02   & 8.2526e-02     &    $\times$    \\
DeepONet & 100 &  101 &   2.8133e-05  &  1.1981e-01  &  4.5926e+02 (20 paths)\\

EMP       & 1000 & 101 & 5.7222e-02    &  1.0523e+00     &    $\times$    \\
DeepONet & 1000 &  101 & 3.2926e-04   &   3.9296e-01    &  4.5926e+02 (20 paths)\\

EMP       & 10000 & 101 & 4.8471e-01   &   9.3683e+01  &    $\times$    \\
DeepONet & 10000 &  101 & 4.6968e-05   &   3.5034e+00  &  4.5926e+02 (20 paths) \\

\hline
\end{tabular}%
}
\caption{Comparison of Errors}
\end{table}
It is worth noting that as $M$ and $N$ increase, the computation time of the EPM method increases at a faster rate, due to the fact that there are too many points to be computed at each iteration step, but DeepONet method does not have this problem. In addition to this, DeepONet can achieve a single Brownian motion orbit to get the corresponding solution directly, which is not possible with EPM.
 
\section{Conclusion}
We propose a novel approach to solve the SDE numerically, i.e. learning the mapping $\bf{F}$ between the solution of SDE and  driven Brownian motion. This mapping can be adapted to various scales and requires minimal training paths to be learned. And we illustrate the multi-scale properties, as well as the strong generalisation ability of the model, through different experimental examples. Finally we use this method for the  McKean-Vlasov SDE to illustrate the application of this method in fast sampling.
\\





\bibliographystyle{elsarticle-harv} 
\bibliography{reference}

\begin{thebibliography}{32}
\expandafter\ifx\csname natexlab\endcsname\relax\def\natexlab#1{#1}\fi
\providecommand{\url}[1]{\texttt{#1}}
\providecommand{\href}[2]{#2}
\providecommand{\path}[1]{#1}
\providecommand{\DOIprefix}{doi:}
\providecommand{\ArXivprefix}{arXiv:}
\providecommand{\URLprefix}{URL: }
\providecommand{\Pubmedprefix}{pmid:}
\providecommand{\doi}[1]{\href{http://dx.doi.org/#1}{\path{#1}}}
\providecommand{\Pubmed}[1]{\href{pmid:#1}{\path{#1}}}
\providecommand{\bibinfo}[2]{#2}
\ifx\xfnm\relax \def\xfnm[#1]{\unskip,\space#1}\fi
\bibitem[{Antonelli and Kohatsu(2002)}]{Antonelli2002RateOC}
\bibinfo{author}{Antonelli, F.}, \bibinfo{author}{Kohatsu, A.}, \bibinfo{year}{2002}.
\newblock \bibinfo{title}{Rate of convergence of a particle method to the solution of the mckean-vlasov's equation}.
\newblock \bibinfo{journal}{Annals of Applied Probability} \bibinfo{volume}{12}, \bibinfo{pages}{423--476}.
\bibitem[{Bao et~al.(2021)Bao, Reisinger, Ren and Stockinger}]{Bao2021FirstorderCO}
\bibinfo{author}{Bao, J.}, \bibinfo{author}{Reisinger, C.}, \bibinfo{author}{Ren, P.}, \bibinfo{author}{Stockinger, W.}, \bibinfo{year}{2021}.
\newblock \bibinfo{title}{First-order convergence of milstein schemes for mckean–vlasov equations and interacting particle systems}.
\newblock \bibinfo{journal}{Proceedings of the Royal Society A} \bibinfo{volume}{477}.
\bibitem[{Bossy and Talay(1997)}]{Bossy1997ASP}
\bibinfo{author}{Bossy, M.}, \bibinfo{author}{Talay, D.}, \bibinfo{year}{1997}.
\newblock \bibinfo{title}{A stochastic particle method for the mckean-vlasov and the burgers equation}.
\newblock \bibinfo{journal}{Math. Comput.} \bibinfo{volume}{66}, \bibinfo{pages}{157--192}.
\bibitem[{Creswell et~al.(2018)Creswell, White, Dumoulin, Arulkumaran, Sengupta and Bharath}]{creswell2018generative}
\bibinfo{author}{Creswell, A.}, \bibinfo{author}{White, T.}, \bibinfo{author}{Dumoulin, V.}, \bibinfo{author}{Arulkumaran, K.}, \bibinfo{author}{Sengupta, B.}, \bibinfo{author}{Bharath, A.A.}, \bibinfo{year}{2018}.
\newblock \bibinfo{title}{Generative adversarial networks: An overview}.
\newblock \bibinfo{journal}{IEEE signal processing magazine} \bibinfo{volume}{35}, \bibinfo{pages}{53--65}.
\bibitem[{Dockhorn et~al.(2021)Dockhorn, Vahdat and Kreis}]{dockhorn2021score}
\bibinfo{author}{Dockhorn, T.}, \bibinfo{author}{Vahdat, A.}, \bibinfo{author}{Kreis, K.}, \bibinfo{year}{2021}.
\newblock \bibinfo{title}{Score-based generative modeling with critically-damped langevin diffusion}.
\newblock \bibinfo{journal}{ArXiv} \bibinfo{volume}{abs/2112.07068}.
\newblock \URLprefix \url{https://api.semanticscholar.org/CorpusID:245131359}.
\bibitem[{Durmus et~al.(2019)Durmus, Majewski and Miasojedow}]{durmus2019analysis}
\bibinfo{author}{Durmus, A.}, \bibinfo{author}{Majewski, S.}, \bibinfo{author}{Miasojedow, B.}, \bibinfo{year}{2019}.
\newblock \bibinfo{title}{Analysis of langevin monte carlo via convex optimization}.
\newblock \bibinfo{journal}{The Journal of Machine Learning Research} \bibinfo{volume}{20}, \bibinfo{pages}{2666--2711}.
\bibitem[{Garg et~al.(2022)Garg, Gupta and Chakraborty}]{Garg2022AssessmentOD}
\bibinfo{author}{Garg, S.}, \bibinfo{author}{Gupta, H.}, \bibinfo{author}{Chakraborty, S.L.}, \bibinfo{year}{2022}.
\newblock \bibinfo{title}{Assessment of deeponet for reliability analysis of stochastic nonlinear dynamical systems}.
\newblock \bibinfo{journal}{ArXiv} \bibinfo{volume}{abs/2201.13145}.
\newblock \URLprefix \url{https://api.semanticscholar.org/CorpusID:246431158}.
\bibitem[{Girolami and Calderhead(2011)}]{girolami2011riemann}
\bibinfo{author}{Girolami, M.}, \bibinfo{author}{Calderhead, B.}, \bibinfo{year}{2011}.
\newblock \bibinfo{title}{Riemann manifold langevin and hamiltonian monte carlo methods}.
\newblock \bibinfo{journal}{Journal of the Royal Statistical Society Series B: Statistical Methodology} \bibinfo{volume}{73}, \bibinfo{pages}{123--214}.
\bibitem[{Goswami et~al.(2022)Goswami, Kontolati, Shields and Karniadakis}]{Goswami2022DeepTL}
\bibinfo{author}{Goswami, S.}, \bibinfo{author}{Kontolati, K.}, \bibinfo{author}{Shields, M.D.}, \bibinfo{author}{Karniadakis, G.E.}, \bibinfo{year}{2022}.
\newblock \bibinfo{title}{Deep transfer learning for partial differential equations under conditional shift with deeponet}.
\newblock \bibinfo{journal}{ArXiv} \bibinfo{volume}{abs/2204.09810}.
\newblock \URLprefix \url{https://api.semanticscholar.org/CorpusID:248299981}.
\bibitem[{Grube(2023)}]{grube2023strong}
\bibinfo{author}{Grube, S.}, \bibinfo{year}{2023}.
\newblock \bibinfo{title}{Strong solutions to mckean--vlasov sdes with coefficients of nemytskii-type}.
\newblock \bibinfo{journal}{Electronic Communications in Probability} \bibinfo{volume}{28}, \bibinfo{pages}{1--13}.
\bibitem[{Hochreiter and Schmidhuber(1997)}]{Hochreiter1997LongSM}
\bibinfo{author}{Hochreiter, S.}, \bibinfo{author}{Schmidhuber, J.}, \bibinfo{year}{1997}.
\newblock \bibinfo{title}{Long short-term memory}.
\newblock \bibinfo{journal}{Neural Computation} \bibinfo{volume}{9}, \bibinfo{pages}{1735--1780}.
\newblock \URLprefix \url{https://api.semanticscholar.org/CorpusID:1915014}.
\bibitem[{Ikeda and Watanabe(2014)}]{ikeda2014stochastic}
\bibinfo{author}{Ikeda, N.}, \bibinfo{author}{Watanabe, S.}, \bibinfo{year}{2014}.
\newblock \bibinfo{title}{Stochastic differential equations and diffusion processes}.
\newblock \bibinfo{publisher}{Elsevier}.
\bibitem[{Kallenberg(2021)}]{Kallenberg2021FoundationsOM}
\bibinfo{author}{Kallenberg, O.}, \bibinfo{year}{2021}.
\newblock \bibinfo{title}{Foundations of modern probability}.
\newblock \bibinfo{journal}{Probability Theory and Stochastic Modelling} \URLprefix \url{https://api.semanticscholar.org/CorpusID:120800745}.
\bibitem[{Kloeden et~al.(1992)Kloeden, Platen, Kloeden and Platen}]{kloeden1992stochastic}
\bibinfo{author}{Kloeden, P.E.}, \bibinfo{author}{Platen, E.}, \bibinfo{author}{Kloeden, P.E.}, \bibinfo{author}{Platen, E.}, \bibinfo{year}{1992}.
\newblock \bibinfo{title}{Stochastic differential equations}.
\newblock \bibinfo{publisher}{Springer}.
\bibitem[{Lanthaler et~al.(2021)Lanthaler, Mishra and Karniadakis}]{Lanthaler2021ErrorEF}
\bibinfo{author}{Lanthaler, S.}, \bibinfo{author}{Mishra, S.}, \bibinfo{author}{Karniadakis, G.E.}, \bibinfo{year}{2021}.
\newblock \bibinfo{title}{Error estimates for deeponets: A deep learning framework in infinite dimensions}.
\newblock \bibinfo{journal}{ArXiv} \bibinfo{volume}{abs/2102.09618}.
\newblock \URLprefix \url{https://api.semanticscholar.org/CorpusID:231979185}.
\bibitem[{Lin et~al.(2020)Lin, Li, Lu, Cai, Maxey and Karniadakis}]{Lin2020OperatorLF}
\bibinfo{author}{Lin, C.}, \bibinfo{author}{Li, Z.}, \bibinfo{author}{Lu, L.}, \bibinfo{author}{Cai, S.}, \bibinfo{author}{Maxey, M.R.}, \bibinfo{author}{Karniadakis, G.E.}, \bibinfo{year}{2020}.
\newblock \bibinfo{title}{Operator learning for predicting multiscale bubble growth dynamics.}
\newblock \bibinfo{journal}{The Journal of chemical physics} \bibinfo{volume}{154 10}, \bibinfo{pages}{104118}.
\newblock \URLprefix \url{https://api.semanticscholar.org/CorpusID:229363556}.
\bibitem[{Liu(2022)}]{Liu2022ParticleMA}
\bibinfo{author}{Liu, Y.}, \bibinfo{year}{2022}.
\newblock \bibinfo{title}{Particle method and quantization-based schemes for the simulation of the mckean-vlasov equation}.
\newblock \bibinfo{journal}{ArXiv} \bibinfo{volume}{abs/2212.14853}.
\bibitem[{Lu et~al.(2019)Lu, Jin and Karniadakis}]{lu2019DeepONet}
\bibinfo{author}{Lu, L.}, \bibinfo{author}{Jin, P.}, \bibinfo{author}{Karniadakis, G.E.}, \bibinfo{year}{2019}.
\newblock \bibinfo{title}{Deeponet: Learning nonlinear operators for identifying differential equations based on the universal approximation theorem of operators}.
\newblock \bibinfo{journal}{arXiv preprint arXiv:1910.03193} .
\bibitem[{Mao(2007)}]{mao2007stochastic}
\bibinfo{author}{Mao, X.}, \bibinfo{year}{2007}.
\newblock \bibinfo{title}{Stochastic differential equations and applications}.
\newblock \bibinfo{publisher}{Elsevier}.
\bibitem[{Oksendal(2013)}]{oksendal2013stochastic}
\bibinfo{author}{Oksendal, B.}, \bibinfo{year}{2013}.
\newblock \bibinfo{title}{Stochastic differential equations: an introduction with applications}.
\newblock \bibinfo{publisher}{Springer Science \& Business Media}.
\bibitem[{Reisinger and Tsianni(2023)}]{reisinger2023convergence}
\bibinfo{author}{Reisinger, C.}, \bibinfo{author}{Tsianni, M.O.}, \bibinfo{year}{2023}.
\newblock \bibinfo{title}{Convergence of the euler--maruyama particle scheme for a regularised mckean--vlasov equation arising from the calibration of local-stochastic volatility models}.
\newblock \bibinfo{journal}{arXiv preprint arXiv:2302.00434} .
\bibitem[{Rezende and Mohamed(2015)}]{rezende2015variational}
\bibinfo{author}{Rezende, D.}, \bibinfo{author}{Mohamed, S.}, \bibinfo{year}{2015}.
\newblock \bibinfo{title}{Variational inference with normalizing flows}, in: \bibinfo{booktitle}{International conference on machine learning}, \bibinfo{organization}{PMLR}. pp. \bibinfo{pages}{1530--1538}.
\bibitem[{Schmuland et~al.(2008)Schmuland, Zhang and Edmonton}]{Schmuland2008YamadaWatanabeTF}
\bibinfo{author}{Schmuland, B.}, \bibinfo{author}{Zhang, X.}, \bibinfo{author}{Edmonton, C.}, \bibinfo{year}{2008}.
\newblock \bibinfo{title}{Yamada-watanabe theorem for stochastic evolution equations in infinite dimensions}.
\newblock \bibinfo{journal}{Condensed Matter Physics} \bibinfo{volume}{11}, \bibinfo{pages}{247}.
\bibitem[{Song et~al.(2020)Song, Sohl-Dickstein, Kingma, Kumar, Ermon and Poole}]{song2020score}
\bibinfo{author}{Song, Y.}, \bibinfo{author}{Sohl-Dickstein, J.N.}, \bibinfo{author}{Kingma, D.P.}, \bibinfo{author}{Kumar, A.}, \bibinfo{author}{Ermon, S.}, \bibinfo{author}{Poole, B.}, \bibinfo{year}{2020}.
\newblock \bibinfo{title}{Score-based generative modeling through stochastic differential equations}.
\newblock \bibinfo{journal}{ArXiv} \bibinfo{volume}{abs/2011.13456}.
\newblock \URLprefix \url{https://api.semanticscholar.org/CorpusID:227209335}.
\bibitem[{Tan and Chen(2022)}]{Tan2022EnhancedDF}
\bibinfo{author}{Tan, L.}, \bibinfo{author}{Chen, L.}, \bibinfo{year}{2022}.
\newblock \bibinfo{title}{Enhanced deeponet for modeling partial differential operators considering multiple input functions}.
\newblock \bibinfo{journal}{ArXiv} \bibinfo{volume}{abs/2202.08942}.
\newblock \URLprefix \url{https://api.semanticscholar.org/CorpusID:246996801}.
\bibitem[{Tappe(2013)}]{Tappe2013TheYT}
\bibinfo{author}{Tappe, S.}, \bibinfo{year}{2013}.
\newblock \bibinfo{title}{The yamada-watanabe theorem for mild solutions to stochastic partial differential equations}.
\newblock \bibinfo{journal}{Electronic Communications in Probability} \bibinfo{volume}{18}, \bibinfo{pages}{1--13}.
\newblock \URLprefix \url{https://api.semanticscholar.org/CorpusID:18398451}.
\bibitem[{Wang et~al.(2021)Wang, Wang and Perdikaris}]{Wang2021LearningTS}
\bibinfo{author}{Wang, S.}, \bibinfo{author}{Wang, H.}, \bibinfo{author}{Perdikaris, P.}, \bibinfo{year}{2021}.
\newblock \bibinfo{title}{Learning the solution operator of parametric partial differential equations with physics-informed deeponets}.
\newblock \bibinfo{journal}{Science Advances} \bibinfo{volume}{7}.
\newblock \URLprefix \url{https://api.semanticscholar.org/CorpusID:232307280}.
\bibitem[{Yang et~al.(2022)Yang, Zhang, Hong, Xu, Zhao, Shao, Zhang, Yang and Cui}]{Yang2022DiffusionMA}
\bibinfo{author}{Yang, L.}, \bibinfo{author}{Zhang, Z.}, \bibinfo{author}{Hong, S.}, \bibinfo{author}{Xu, R.}, \bibinfo{author}{Zhao, Y.}, \bibinfo{author}{Shao, Y.}, \bibinfo{author}{Zhang, W.}, \bibinfo{author}{Yang, M.H.}, \bibinfo{author}{Cui, B.}, \bibinfo{year}{2022}.
\newblock \bibinfo{title}{Diffusion models: A comprehensive survey of methods and applications}.
\newblock \bibinfo{journal}{ArXiv} \bibinfo{volume}{abs/2209.00796}.
\bibitem[{Zhang et~al.(2022a)Zhang, Zhang and Lin}]{zhang2022multiauto}
\bibinfo{author}{Zhang, J.}, \bibinfo{author}{Zhang, S.}, \bibinfo{author}{Lin, G.}, \bibinfo{year}{2022}a.
\newblock \bibinfo{title}{Multiauto-deeponet: A multi-resolution autoencoder deeponet for nonlinear dimension reduction, uncertainty quantification and operator learning of forward and inverse stochastic problems}.
\newblock \bibinfo{journal}{arXiv preprint arXiv:2204.03193} .
\bibitem[{Zhang et~al.(2022b)Zhang, Zhang and Lin}]{Zhang2022MultiAutoDeepONetAM}
\bibinfo{author}{Zhang, J.}, \bibinfo{author}{Zhang, S.}, \bibinfo{author}{Lin, G.}, \bibinfo{year}{2022}b.
\newblock \bibinfo{title}{Multiauto-deeponet: A multi-resolution autoencoder deeponet for nonlinear dimension reduction, uncertainty quantification and operator learning of forward and inverse stochastic problems}.
\newblock \bibinfo{journal}{ArXiv} \bibinfo{volume}{abs/2204.03193}.
\newblock \URLprefix \url{https://api.semanticscholar.org/CorpusID:248006371}.
\bibitem[{Zhang and Chen(2021)}]{zhang2021diffusion}
\bibinfo{author}{Zhang, Q.}, \bibinfo{author}{Chen, Y.}, \bibinfo{year}{2021}.
\newblock \bibinfo{title}{Diffusion normalizing flow}.
\newblock \bibinfo{journal}{Advances in Neural Information Processing Systems} \bibinfo{volume}{34}, \bibinfo{pages}{16280--16291}.
\bibitem[{Zhao(2014)}]{Zhao2014YamadaWatanabeTF}
\bibinfo{author}{Zhao, H.}, \bibinfo{year}{2014}.
\newblock \bibinfo{title}{Yamada-watanabe theorem for stochastic evolution equation driven by poisson random measure}.
\newblock \bibinfo{journal}{International Scholarly Research Notices} \bibinfo{volume}{2014}, \bibinfo{pages}{1--7}.
\newblock \URLprefix \url{https://api.semanticscholar.org/CorpusID:71103122}.

\end{thebibliography}



\end{document}